\documentclass{amsart}

\usepackage[T1]{fontenc}	% see http://tex.stackexchange.com/questions/664/why-should-i-use-usepackaget1fontenc
\usepackage{lmodern}		% see http://tex.stackexchange.com/questions/1291/why-are-bitmap-fonts-used-automatically 
\usepackage[utf8]{inputenc}	% see http://tex.stackexchange.com/questions/44694/fontenc-vs-inputenc
\usepackage[english]{babel}	% for better hyphenation of english words

\usepackage{amsfonts,amssymb,amsmath,tikz}

\usepackage{xcolor}

\usepackage{mleftright}\mleftright  % (sometimes) fixes spacing of \left and \right

\usepackage[bookmarks=true, bookmarksnumbered=true, bookmarksopen=true,
	unicode=false, pdftoolbar=true, pdfmenubar=true, pdffitwindow=false, pdfstartview={FitH},
	pdftitle={}, pdfauthor={}, pdfsubject={}, pdfcreator={}, pdfproducer={}, pdfkeywords={},
	pdfnewwindow=true, colorlinks=true, linkcolor=black, citecolor=black, filecolor=black, urlcolor=black]{hyperref}

\newcommand{\R}{\mathbb{R}}
\newcommand{\N}{\mathbb{N}}
\newcommand{\E}{\mathbb{E}}
\newcommand{\Pro}{\mathbb{P}}
\newcommand{\Var}{\operatorname{Var}}

\renewcommand{\S}{\ensuremath{{\mathbb{S}}}}
\newcommand{\KS}{\operatorname{dist_{K}}}

\DeclareMathOperator{\conv}{conv}
\DeclareMathOperator{\vol}{vol}

\usepackage{amsthm}

\newtheorem{thm}{Theorem}[section]
\newtheorem{cor}[thm]{Corollary}
\newtheorem{lemma}[thm]{Lemma}

\newtheorem{rmk}[thm]{Remark}

\begin{document}

%%%%%%%%%%%%%%%%%%%%%%%%%%%%%%%%%%%%%%%%%%%%%5

\title{Asymptotic normality for random simplices and convex bodies in high dimensions}

\author[D.~Alonso-Guti\'errez]{D.\ Alonso-Guti\'errez}
\address{University of Zaragoza, Spain} 
\email{alonsod@unizar.es}

\author[F.~Besau]{F.\ Besau}
\address{Vienna University of Technology, Austria} 
\email{florian.besau@tuwien.ac.at}

\author[J.~Grote]{J.\ Grote}
\address{University of Ulm, Germany} 
\email{julian.grote@uni-ulm.de}

\author[Z.~Kabluchko]{Z.\ Kabluchko}
\address{University of M\"unster, Germany} 
\email{zakhar.kabluchko@wwu.de}

\author[M.~Reitzner]{M.\ Reitzner}
\address{University of Osnabrück, Germany} 
\email{matthias.reitzner@uni-osnabrueck.de}

\author[C.~Th\"ale]{C.\ Th\"ale}
\address{Ruhr University Bochum, Germany} 
\email{christoph.thaele@rub.de}

\author[B.-H.~Vritsiou]{B.-H.\ Vritsiou}
\address{University of Alberta in Edmonton, Canada} 
\email{vritsiou@ualberta.ca}

\author[E.~Werner]{E.\ Werner}
\address{Case Western Reserve University, USA} 
\email{elisabeth.werner@case.edu}


\date{}

\begin{abstract}
Central limit theorems for the log-volume of a class of random convex bodies in $\mathbb{R}^n$ are obtained in the high-dimensional regime, that is, as $n\to\infty$. In particular, the case of random simplices pinned at the origin and simplices where all vertices are generated at random is investigated. 
The coordinates of the generating vectors are assumed to be independent and identically distributed with subexponential tails. 
In addition, asymptotic normality is established also for random convex bodies (including random simplices pinned at the origin) when the spanning vectors are distributed according to a radially symmetric probability measure on the $n$-dimensional $\ell_p$-ball. In particular, this includes the cone and the uniform probability measure.
\end{abstract}

\maketitle

\section{Introduction and main results}

\subsection{Motivation}

Central limit theorems for random polytopes in $\R^n$ are widely known if the space dimension $n$ is kept fixed, while the number of generating points tends to infinity. We refer, for example, to the survey articles of B\'ar\'any \cite{BaranySurvey}, Hug \cite{HugSurvey} and Reitzner \cite{ReitznerSurvey} for results in this direction and for further references. In the present paper we investigate the case where the number of generating points is essentially equal to the space dimension and both tend to infinity \textit{simultaneously}. To be more precise, we consider the case of random $n$-dimensional simplices in $\R^n$, where we distinguish between the case of $(n+1)$ generating points chosen at random, or the case where we only have $n$ random points and the $(n+1)$st vertex is fixed at the origin. The latter construction is called a \emph{pinned simplex} in the following. 
Asymptotic normality for the log-volume of random simplices in high dimensions has previously been considered by Ruben \cite{Ruben1977}, Maehara \cite{Maehara1980} and Mathai \cite{Mathai82}. Note however, that in their results the dimension of the random (pinned) simplices is kept fixed, while the space dimension $n$ tends to infinity. For Gaussian and so-called beta simplices Eichelsbacher and Knichel \cite{EichelsbacherKnichel} and Grote, Kabluchko and Th\"ale \cite{GroteKabluchkoThaele} recently also studied a number of probabilistic limit theorems where the simplex dimension tends to infinity as well.

Our main result is a central limit theorem for the log-volume of a random $n$-dimensional simplix in $\R^n$, see Theorem \ref{thm:CLTProductMeasurePara} below. 
For an $n$-dimensional random simplex that is pinned, it is known that its volume is determined by the absolute value of the determinant of the matrix whose columns are given by the generating vectors from the origin. As a consequence, if these columns are filled by independent and identically distributed (i.i.d.) random variables, a central limit theorem for the log-volume of random pinned simplex follows from the central limit theorem for random determinants with i.i.d.\ entries established by Nguyen and Vu \cite{Nguyen-Vu-2014}.

The same arguments cannot directly be applied if the coordinates of the generating points are not independent, that is, for example if the points are chosen with respect to a probability measure in the $\ell_p$-ball $B_p^n\subset \R^n$ with $p\neq\infty$. However, we still succeed in establishing a central limit theorem for random pinned simplices in the $\ell_p$-ball for certain radially symmetric probability measures, which include in particular the uniform probability measure and the cone-volume measure, see Theorem \ref{thm:CLTlpCASE} below. In our proof we employ different tools, most notably a Schechtman-Zinn-type probabilistic representation of Barthe, Gu\'edon, Mendelson and Naor \cite{BartheEtal2005}, which allow us to relate the log-volume of the random pinned simplex to the determinant of a matrix with independent entries.
Hence, although the coordinates of the generating vectors are now no longer independent, at the core of our argument we can still rely on the central limit theorem for the determinant of random matrices with i.i.d.\ entries.

\subsection{Main results: the case of independent coordinates}

Let $\mu$ be a probability measure on $\R^n$ ($n\geq 1$) and let $X_0,\dotsc,X_n$ be independent random vectors distributed according to $\mu$. We define the random simplex
\begin{equation*}
	\Sigma_n   := \conv\bigl( \{X_0,X_1, \dotsc, X_n\} \bigr),
\end{equation*}
as well as the random pinned simplex
\begin{equation*}
	\Sigma_n^0 := \conv\bigl( \{0,X_1,\dotsc, X_n\}\bigr).
\end{equation*}
In what follows, we shall focus on the case where the coordinates $\xi^i_j$, $j\in\{1,\dotsc,n\}$, of the random vectors $X_i=(\xi^i_1,\ldots,\xi^i_n)$ are independent copies of a random variable $\xi$. Furthermore, we assume that the random variable $\xi$ is symmetric, has variance one and subexponential tails with exponent $\alpha>0$. By the latter we mean that there are constants $c_1,c_2>0$ such that
\begin{equation}\label{eq:AssumptionTails}
	\Pro\bigl(|\xi|\geq t^\alpha \bigr)\leq c_1 e^{-c_2t}, \qquad t>0.
\end{equation}
Examples are the uniform distribution on the cube $[-\sqrt{3},\sqrt{3}]^n$, the uniform distribution on the discrete cube $\{-1,+1\}^n$,
the two-sided exponential distribution, standard Gaussian distribution on $\R^n$ or, more generally,
the $p$-generalized Gaussian distribution with density proportional to
$e^{-|t|^p/a}$ (for an appropriate choice of $a>0$) for any $p>0$.

In our first result we establish a central limit theorem for the log-volume of the high-dimensional random simplices $\Sigma_n$ and $\Sigma_n^0$, as $n\to\infty$. In the following $Z\sim\mathcal{N}(0,1)$ always denotes a standard Gaussian random variable and $\overset{d}{\longrightarrow}$ indicates convergence in distribution. 

\begin{thm}[CLT for random simplices]\label{thm:CLTProductMeasurePara}
Let $\xi$ be a symmetric random variable with variance one and subexponential tails with exponent $\alpha>0$.

\begin{enumerate}
		\item[i)] Assume $\Sigma_n   := \conv\bigl( \{X_0,X_1, \dotsc, X_n\} \bigr)$ is a random simplex in $\R^n$ with $X_i$ having i.i.d. coordinates $\xi_j^{i} \sim \xi$. Then
	\begin{equation*}
		\frac{\ln \vol_n(\Sigma_n) +\frac{n}{2}\ln n - \frac{n}{2}+\frac{1}{4} \ln n}{\sqrt{\frac{1}{2}\ln n}} \overset{d}{\longrightarrow} Z, \qquad \text{as $n\to\infty$.}
	\end{equation*}
\item[ii)] Assume $\Sigma_n^0  := \conv\bigl( \{0,X_1, \dotsc, X_n\} \bigr)$ is a random simplex in $\R^n$ with $X_i$ having i.i.d. coordinates $\xi_j^{i} \sim \xi$. Then
	\begin{equation*}
		\frac{\ln \vol_n(\Sigma_n^0) + \frac{n}{2} \ln  n -\frac{n}{2}+\frac{3}{4} \ln n}{\sqrt{\frac{1}{2}\ln n}} \overset{d}{\longrightarrow} Z, \qquad \text{as $n\to\infty$.}
	\end{equation*}
\end{enumerate}
\end{thm}

Part ii) of Theorem \ref{thm:CLTProductMeasurePara} can be reformulated for the parallelotope spanned by the vertices of the pinned simplex from the origin. Even more generally, for a convex body $K\subset \R^n$, that is a compact convex subset with non-empty interior, and $n$ random points $X_1,\dotsc, X_n \in \R^n$ we define the random convex body
\begin{equation}\label{eqn:randConst}
	\Xi_n(K) := K[X_1,\dotsc,X_n] = \left\{ \sum_{i=1}^n y_i X_i : (y_1,\dotsc,y_n)\in K\right\}.
\end{equation}
We note that for each $n\in\N$ and each convex body $K\subset\R^n$, $\Xi_n(K)$ is a random closed set in the usual sense of stochastic geometry, cf.\ \cite[Chapter 16]{Kallenberg}. In particular, this implies that that the volume $\vol_n(\Xi_n(K))$ of $\Xi_n(K)$ is an ordinary random variable. As observed by Paouris and Pivovarov \cite{ PaourisPivovarov:2012, PaourisPivovarov:2013}, this concept generalizes a number of common constructions. Namely,
\begin{enumerate}
	\item[a)] if $K$ is the standard simplex 
		\begin{equation*}
			T^n = \left\{ (x_1,\dotsc,x_n) \in \R^n: x_i\geq 0\text{ and } \sum_{i=1}^n x_i \leq 1\right\},
		\end{equation*}
		then $\Xi_n(T^n)$ coincides with the pinned simplex $\Sigma_n^0$.
	\item[b)] if $K$ is the unit cube $C^n=[0,1]^n$, then $\Xi_n(C^n)$ is the parallelotope spanned by $(X_i)_{i=1}^n$ from the origin.
	\item[c)] if $K = B_\infty^n = [-1,1]^n$ is the symmetric cube, then $\Xi_n(B_\infty^n)$ is the zonotope generated by the segments $[-X_i,X_i]$, i.e.,
	\begin{equation*}
		\Xi_n(B_\infty^n)  = \left\{\sum_{i=1}^n \lambda_i X_i : \lambda_i\in [-1,1]\right\}.
	\end{equation*}
	\item[d)] if $K = B_1^n$ is the cross-polytope, then $\Xi_n(B^n_1)$ is the symmetric convex hull of the $2n$ points $\{\pm X_i: i=1,\dotsc,n\}$.
	\item[e)] if $K = B_2^n$ is the unit ball, then $\Xi_n(B^n_2)$ is an ellipsoid, that is, it is the image of the unit ball under the linear map whose matrix is generated by the random points $(X_i)_{i=1}^n$.
\end{enumerate}

As a generalization of part ii) of Theorem \ref{thm:CLTProductMeasurePara} we obtain the following central limit theorem for the log-volume of the random convex bodies $\Xi_n(K)$. In the following we denote by $\KS(-,-)$ the Kolmogorov distance between random variables, that is, for two random variables $X,Y$ we have
\begin{equation*}
	\KS(X,Y) := \sup_{t\in \R}\, \bigl| \Pro(X\leq t)-\Pro(Y\leq t)\bigr|.
\end{equation*}
Note that convergence in the Kolmogorov distance implies convergence in distribution. Also, by $o(1)$ we denote some sequence $(a_n)_{n\in\mathbb{N}}$ with $a_n\to 0$, as $n\to \infty$.

\pagebreak
\begin{thm}\label{thm:CLTProductMeasureGen}
	Let $\xi$ be a symmetric random variable with variance one and subexponential tails with exponent $\alpha>0$. Let $(K_n)_{n\in\mathbb{N}}$ be a sequence of convex bodies such that $K_n\subset \R^n$.
	If $X_1,\dotsc,X_n$ are random points in $\R^n$ with i.i.d. coordinates $\xi_j^{i} \sim \xi$, and $\Xi_n(K_n)$ is the random convex body as defined by \eqref{eqn:randConst}, then
	\begin{equation*}
		S_n:=\frac{\ln \vol_n(\Xi_n(K_n)) - \ln \vol_n(K_n) - \frac{n}{2} \ln n + \frac{n}{2} + \frac{1}{4} \ln n}{\sqrt{\frac{1}{2}\ln n}} \overset{d}{\longrightarrow} Z, \qquad \text{as $n\to\infty$.}
	\end{equation*}
	More precisely, we have that
	\begin{equation*}
		\KS(S_n,Z)\leq (\ln n)^{-\frac{1}{3} + o(1)}.
	\end{equation*}
\end{thm}

As an application of Theorem \ref{thm:CLTProductMeasureGen} we may  revisit the special cases a) -- d) of the random convex bodies $\Xi_n(K_n)$ mentioned above. The resulting central limit theorems are summarized in Table \ref{table:CLT}. In particular, taking $K_n=T^n$ for each $n\in\N$, we also obtain part ii) of Theorem \ref{thm:CLTProductMeasurePara}.

\begin{table}[t]
\begin{tabular}{|c|c||c|}
\hline
$\parbox[0pt][2em][c]{0cm}{}K_n$ & $\vol_n(K_n)$ & central limit theorem\\
\hline
\hline
$\parbox[0pt][3em][c]{0cm}{}T^n$ & $\frac{1}{n!}$ & $\frac{\ln \vol_n(\Xi_n(T^n)) + \frac{n}{2} \ln  n - \frac{n}{2} +\frac{3}{4} \ln n}{\sqrt{\frac{1}{2}\ln n}} \overset{d}{\longrightarrow} Z$\\
\hline
$\parbox[0pt][3em][c]{0cm}{}C^n$ & $1$ & $\frac{\ln \vol_n(\Xi_n(C^n)) - \frac{n}{2} \ln n + \frac{n}{2} +\frac{1}{4} \ln n}{\sqrt{\frac{1}{2}\ln n}}\overset{d}{\longrightarrow} Z$\\
\hline
$\parbox[0pt][3em][c]{0cm}{}B_\infty^n$ & $2^n$ & $\frac{\ln \vol_n(\Xi_n(B_\infty^n)) - \frac{n}{2} \ln n -(\ln 2- \frac{1}{2}) n + \frac{1}{4} \ln n}{\sqrt{\frac{1}{2}\ln n}}\overset{d}{\longrightarrow} Z$\\
\hline
$\parbox[0pt][3em][c]{0cm}{}B_1^\infty$ & $\frac{2^n}{n!}$ & $\frac{\ln \vol_n(\Xi_n(B_1^n)) + \frac{n}{2} \ln - (\ln 2 + \frac{1}{2}) n + \frac{3}{4} \ln n}{\sqrt{\frac{1}{2}\ln n}}\overset{d}{\longrightarrow} Z$\\
\hline
\end{tabular}
\vspace{2mm}
\caption{Special cases of the central limit theorem for the log-volume of random covex bodies (Theorem \ref{thm:CLTProductMeasureGen}). Here, $Z$ is a standard Gaussian random variable.}\label{table:CLT}
\end{table}

\subsection{Main results: the case of \texorpdfstring{$\ell_p$}{Lp}-balls}

For $0<p\leq \infty$ the $n$-dimensional $\ell_p$-ball $B_p^n\subset \R^n$ is defined as
\begin{equation*}
	B_p^n := \left\{ x\in\R^n:\|x\|_p\leq 1 \right\},
\end{equation*}
where the $p$-norm (or quasi-norm if $0<p<1$) of $x=(x_1,\dotsc,x_n)\in\R^n$ is
\begin{equation*}
	\|x\|_p := \begin{cases}
		\left(\sum\limits_{i=1}^n|x_i|^p\right)^{1/p} &\text{if $0<p<\infty$,}\\
		\max\{|x_1|,\ldots,|x_n|\} &\text{if $p=\infty$.}
	\end{cases}
\end{equation*}
In our next result we consider pinned simplices, denoted by $\Sigma_n^0(\nu)$, which are spanned by the origin and $n$ points $X_1,\ldots,X_n$ chosen at random with respect to a radially symmetric probability measure $\nu=\nu_n(m_n,p)$ on the $n$-dimensional $\ell_p$-ball $B_p^n$. More specifically, $\nu$ belongs to a family of measures including the cone probability measure and the uniform probability measure on $B_p^n$ which is driven by a parameter $m_n \geq 0$. This model contains a number of special cases that are of particular interest (see Theorem 1, Theorem 3, Corollary 3 and Corollary 4 in \cite{BartheEtal2005} as well as the discussion before Theorem 1.1 in \cite{APT17CLT}). Namely, if $m_n=0$, then the random points $X_1,\ldots,X_n$ are distributed according to the cone probability measure on the boundary of $B_p^n$, i.e., the $\ell_p$-sphere in $\R^n$. It is well known that this measure coincides with the normalized surface measure precisely if $p\in\{1,2,\infty\}$. Next, if $m_n=1$, then $X_1,\ldots,X_n$ are selected according to the uniform distribution on $B_p^n$. Finally, if $m_n=m/p$ for some $m\in\N$, then the distribution corresponds to the image of the cone probability measure on $B_p^{n+m}$ under the orthogonal projection onto the first $n$ coordinates. Similarly, if $m_n=1+m/p$, then the points are sampled according to the image of the uniform distribution on $B_p^{n+m}$ under the same projection. We refer to Section 3 for the precise construction of $\nu$. 

\begin{thm}[CLT for random convex bodies in the $\ell_p$-ball]\label{thm:CLTlpCASE}
Let $X_1,\dotsc, X_n$ be $n$ independent random points in the $\ell_p$-ball $B_p^n$ with respect to a probability measure $\nu=\nu_n(m_n,p)$ as defined in Section \ref{sec:nu}.
Let $(K_n)_{n\in\mathbb N}$ be a sequence of convex bodies such that $K_n\subset \R^n$ and $\Xi_n(K_n)$ be the random convex body generated by $\nu$-distributed random points $X_1,\dotsc,X_n$ as defined by \eqref{eqn:randConst}. Then
\begin{equation*}
	S_n:=\frac{\ln \vol_n(\Xi_n(K_n)) - \ln\vol_n(K_n) - \frac{1}{2} \ln\, (n-1)! + \frac{n}{p} \ln(a(m_n+\frac{n}{p}))}{\sqrt{\frac{1}{2}\ln n}} \overset{d}{\longrightarrow} Z,
\end{equation*}
as $n\to\infty$, where $a=\left({\Gamma(1/p)}/{\Gamma(3/p)}\right)^{p/2}$.
More precisely, we have that
\begin{equation*}
	\KS(S_n,Z)\leq (\ln n)^{-\frac{1}{3}+o(1)}.
\end{equation*}
\end{thm}

By choosing for each $n\in\N$, $K_n$ as the standard simplex $T^n$ we obtain the following central limit theorem as a direct corollary to Theorem \ref{thm:CLTlpCASE}.

\begin{cor}[CLT for random pinned simplices in the $\ell_p$-ball]
	Let $\Sigma_n^0(\nu)$ be the random pinned simplex that is spanned by the origin and $n$ independent random points in the $\ell_p$-ball $B_p^n$ which are distributed according to a probability measure $\nu=\nu_n(m_n,p)$ as defined in Section \ref{sec:nu}. Then
	\begin{equation*}
		\frac{\ln \vol_n(\Sigma_n^0(\nu)) +\frac{n}{2} \ln n-\frac{n}{2} +\frac{3}{4} \ln n + \frac{n}{p} \ln(a(m_n+\frac{n}{p}))}{\sqrt{\frac{1}{2}\ln n}} \overset{d}{\longrightarrow} Z, \qquad \text{as $n\to\infty$,}
	\end{equation*}
	where $a=\left({\Gamma(1/p)}/{\Gamma(3/p)}\right)^{p/2}$.
\end{cor}

\subsection{Plan of the paper}

In the next Section we outline the proof of Theorem \ref{thm:CLTProductMeasurePara} and of Theorem \ref{thm:CLTProductMeasureGen}. We will collect the relevant tools along the way. As mentioned in the introduction, the proof of Theorem \ref{thm:CLTProductMeasurePara} essentially relies on the central limit theorem for determinants of random matrices of Nguyen and Vu \cite{Nguyen-Vu-2014} with additional arguments for the non-pinned case. In Section 3 we present the details of the proof of Theorem \ref{thm:CLTlpCASE} and recall the definition of the special measure $\nu=\nu(m_n,p)$. For the proof of Theorem \ref{thm:CLTlpCASE} we especially need the Schechtman-Zinn-type probabilistic representation from \cite{BartheEtal2005}.

\subsection*{Acknowledgment}
The authors started this project within a working group that formed during the Mini-Workshop \textit{Perspectives in High-dimensional Probability and Convexity} at the Mathematisches Forschungsinstitut Oberwolfach (MFO). All support is gratefully acknowledged.

D.A.-G.\ is partially supported by DGA grant \verb|E26_17R| and MINECO grant \verb|MTM2016-77710-P| and IUMA.
F.B.\ was partially supported by the Deutsche For-schungsgemeinschaft (DFG) grant \verb|BE 2484/5-2|.
J.G.\ was supported by DFG via \verb|RTG 2131| ``High-dimensional Phenomena in Probability – Fluctuations and Discontinuity''.
E.W.\ was partially supported by NSF grant \verb|DMS-1811146|.

\section{Proof of Theorem \ref{thm:CLTProductMeasurePara} and Theorem \ref{thm:CLTProductMeasureGen}}
\label{sec:ProofIIDCase}

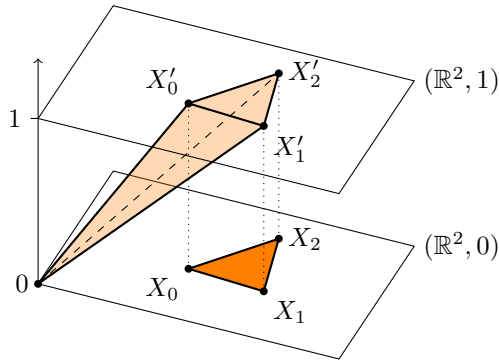
\begin{figure}[t]
	\begin{center}
	\begin{tikzpicture}
	\def\h{2.2} %distance of plane (R^n,1) and (R^n,0);
	
	\draw (0,0) -- (4,-1) -- (5,0.5) -- (1,1.5) -- cycle;
	
	\draw[->] (0,0) -- (0,3);
	\fill (0,0) circle (1.5pt);

	\coordinate (A1) at (2,0.2+\h);
	\coordinate (A2) at (3,-0.1+\h);
	\coordinate (A3) at (3.2,0.6+\h);
	\fill[orange, opacity=0.3] (0,0) -- (A1) -- (A2) -- cycle;
	\fill[orange, opacity=0.3] (A1) -- (A2) -- (A3) -- cycle;
	\draw[thick] (0,0) -- (A1);
	\draw[thick] (0,0) -- (A2);
	\draw[dashed] (0,0) -- (A3);
	\draw[thick] (A1) -- (A2) -- (A3) -- cycle;
	\fill (A1) circle (1.5pt) node[above left] {$X_0'$};
	\fill (A2) circle (1.5pt) node[below right] {$X_1'$};
	\fill (A3) circle (1.5pt) node[right] {$X_2'$};

	\draw (0,0+\h) -- (4,-1+\h) -- (5,0.5+\h) -- (1,1.5+\h) -- cycle;

	\coordinate (B1) at (2,0.2);
	\coordinate (B2) at (3,-0.1);
	\coordinate (B3) at (3.2,0.6);
	\draw[thick, fill=orange] (B1) -- (B2) -- (B3) -- cycle;
	\fill (B1) circle (1.5pt) node[below left] {$X_0$};
	\fill (B2) circle (1.5pt) node[below right] {$X_1$};
	\fill (B3) circle (1.5pt) node[right] {$X_2$};
	
	\draw[dotted] (A1) -- (B1);
	\draw[dotted] (A2) -- (B2);
	\draw[dotted] (A3) -- (B3);

	\draw (0,\h) -- (-0.1,\h) node[left] {$1$};
	\node[left] at (0,0) {$0$};
	\node[right] at (5,0.5+\h) {$(\R^2,1)$};
	\node[right] at (5,0.5) {$(\R^2,0)$};
\end{tikzpicture}
\caption{Illustration of the identification of the random simplex $\Sigma_n$ with a pinned simplex $\Sigma_{n+1}^0$ in the two dimensional case.}\label{fig}
\end{center}
\end{figure}

\subsection{Proof of Theorem \ref{thm:CLTProductMeasureGen} and part ii) of Theorem \ref{thm:CLTProductMeasurePara}}

Let us first recall the central limit theorem of Nguyen and Vu for the log-determinant of random matrices with independent entries.

\begin{thm}[{\cite[Theorem 1.1]{Nguyen-Vu-2014}}]\label{thm:CLTDeterminantNguyenVu}
Let $A_n$ be an $n\times n$ random matrix whose entries are independent random variables with zero mean, variance one and subexponential tails with exponent $\alpha>0$. Further, let $Z\sim\mathcal{N}(0,1)$ be a standard Gaussian random variable. Then,
	\begin{equation*}
		Z_n:=\frac{\ln\,\left|\det A_n\right|-\frac{1}{2}\ln\, (n-1)!}
			{\sqrt{\frac{1}{2}\ln n}} \overset{d}{\longrightarrow} Z,\quad \text{as $n\to \infty$.}
	\end{equation*}
	More precisely, the rate of convergence is
	\begin{equation}\label{eqn:CLTdetStrong}
		\KS(Z_n,Z)
		\leq (\ln n)^{-1/3+o(1)},
	\end{equation}
	for all $n$ large enough.
\end{thm}

Next, we recall the definition \eqref{eqn:randConst} of the random convex bodies $\Xi_n(K_n)$ and observe that
\begin{equation}\label{eqn:randbody_vol}
	\vol_n(\Xi_n(K_n)) = \left|\det(X_1|\dotsc|X_n)\right|\vol_n(K_n),
\end{equation}
see, for example, \cite[Proposition 2.1]{PaourisPivovarov:2013}. This implies that
\begin{align*}
	&\frac{\ln\vol_n(\Xi_n(K_n))-\ln \vol_n(K_n)-\frac{1}{2}\ln\,(n-1)!}{\sqrt{\frac{1}{2}\ln\,n}} \\
	&\qquad\qquad\qquad\qquad =\frac{\ln\, \left|\det(X_1|\dotsc|X_n)\right|-\frac{1}{2}\ln\,(n-1)!}{\sqrt{\frac{1}{2}\ln\,n}}
\end{align*}
and hence Theorem \ref{thm:CLTProductMeasureGen} is a direct consequence of Theorem \ref{thm:CLTDeterminantNguyenVu}.  Moreover, taking $K_n=T^n$ for each $n\in\N$ and recalling that $\Sigma_n^0 = \Xi_n(T^n)$, we may derive part ii) in Theorem \ref{thm:CLTProductMeasurePara} by Stirling's formula,
\begin{equation*}
	\ln n! = n\ln n - n + \frac{1}{2}\ln n + O(1).
\end{equation*}

\subsection{Proof of part i) in Theorem \ref{thm:CLTProductMeasurePara}}
We now prove the central limit theorem for the log-volume of random simplices $\Sigma_n$. First, notice that the volume of $\Sigma_n$ can be identified with the volume of a pinned simplex in $\R^{n+1}$ by the classical projective construction (see also Figure \ref{fig}). Given $n+1$ points $X_0,X_1,\dotsc, X_{n}$ in $\R^n$ we consider the points $X_i' :=(X_i,1) \in \R^{n+1}$ for $i=0,\dotsc,{n}$. We set $\Sigma_{n+1}^0 := \conv(\{0,X_0',\dotsc,X_{n}'\})$ and find that
\begin{equation*}
	\vol_n(\Sigma_n) = (n+1) \vol_{n+1}(\Sigma_{n+1}^0) = \frac{1}{n!} \left|\det(X_0'|\dotsc|X_n')\right|.
\end{equation*}
Now let $Y_i\in \R^{n+1}$ be the vector whose entries are given by the $i$th row of the $(n+1)\times (n+1)$ matrix $(X_0'|\dotsc|X_n')$ for $i=1,\dotsc, n$. Then $Y_i$ is a random vector in $\R^{n+1}$ with independent coordinates distributed like $\xi$. Notice that the last row in the matrix $(X_0'|\dotsc|X_n')$ is just the constant vector $U_{n+1}:=(1,\dotsc,1)\in\R^{n+1}$.
Let $Y_{n+1}$ be another random vector with independent entries distributed like $\xi$. We compare the random matrix $(Y_1|\dotsc|Y_n|U_{n+1})$ with the random matrix $(Y_1|\dotsc|Y_n|Y_{n+1})$. Notice that the latter now has independent and identically distributed entries.
We have
\begin{align*}
	\frac{1}{n!} \left|\det(Y_1|\dotsc|Y_n|U_{n+1})\right| &= \operatorname{dist}(U_{n+1},L_n) \vol_n(\conv \{0,Y_1,\dotsc,Y_n\})\\
	\intertext{and}
	\frac{1}{n!} \left|\det(Y_1|\dotsc|Y_n|Y_{n+1})\right| &= \operatorname{dist}(Y_{n+1},L_n) \vol_n(\conv \{0,Y_1,\dotsc,Y_n\}),
\end{align*}
where $L_n$ is the $n$-dimensional linear subspace spanned by $Y_1,\dotsc,Y_n$ in $\R^{n+1}$ and $\operatorname{dist}(v,L_n)$ denotes the distance of a vector $v\in\R^{n+1}$ to $L_n$.
Collecting all of the above we conclude that
\begin{align*}
	\ln \vol_n(\Sigma_n) &= \ln\, \left|\det(X_0'|\dotsc|X_n')\right| - \ln n! \\
	&= \ln\, \left|\det(Y_1|\dotsc|Y_{n+1})\right| -\ln \operatorname{dist}(Y_{n+1},L_n) + \ln \operatorname{dist}(U_{n+1},L_n)- \ln n!.
\end{align*}
We will show that
\begin{align}\label{eqn:proofAii}
	\frac{\ln \operatorname{dist}(U_{n+1},L_n)}{\sqrt{\frac{1}{2} \ln n}} &\overset{d}{\longrightarrow} 0
	&\text{and} &&
	\frac{\ln \operatorname{dist}(Y_{n+1},L_n)}{\sqrt{\frac{1}{2} \ln n}} &\overset{d}{\longrightarrow} 0,
\end{align}
as $n\to \infty$. Then we may conclude by Slutsky's theorem (see, for example, \cite[Proposition A.42 (b)]{BassBook}) and the central limit theorem for the log-determinant, Theorem \ref{thm:CLTDeterminantNguyenVu}, that
\begin{align*}
	\frac{\ln \vol_n(\Sigma_n) + \frac{1}{2} \ln n!}{\sqrt{\frac{1}{2} \ln n}} 
	&= \frac{\ln\, \left|\det(Y_1|\dotsc|Y_{n+1})\right| - \frac{1}{2} \ln n!}{\sqrt{\frac{1}{2} \ln n}} \\
	&\qquad+ \frac{\ln \operatorname{dist}(U_{n+1},L_n)}{\sqrt{\frac{1}{2} \ln n}} - \frac{\ln \operatorname{dist}(Y_{n+1},L_n)}{\sqrt{\frac{1}{2} \ln n}}
\end{align*}
converges in distribution to a standard Gaussian random variable $Z\sim\mathcal{N}(0,1)$, as $n\to \infty$.

\medskip
To prove \eqref{eqn:proofAii} we need the following two auxiliary results.

\begin{lemma}[{Berry-Esseen inequality \cite[Lemma 8.1]{Nguyen-Vu-2014}}] \label{lem:Berry-Esseen-sup}
	Let $B_n=(b_1,\dotsc,b_n)$ be a random vector whose coordinates are independent copies of a random variable $\xi$ with mean zero, variance one and subexponential tails with exponent $\alpha>0$ and let $V_{n} = (v_1, \ldots, v_n)$ be a fixed unit vector in $\R^{n}$.
	Then there exists a constant $c\in(0,\infty)$ such that
	\begin{equation}\label{eq:BerryEsseenProof2}
		\KS(|\langle V_n, B_n \rangle|,|Z|) \leq c\|V_{n}\|_{\infty}. 
	\end{equation}
\end{lemma}

\begin{lemma}[{\cite[Theorem 1.4]{NguyenVuNormal} for subexponential tails, see Remark \ref{rmk:note} below}]\label{le:NguyenVuNormal}
Suppose $L_n$ is a linear subspace spanned by $n$ independent random vectors $Y_1, \dotsc, Y_n$ in ${\mathbb R}^{n+1}$ each of whose coordinates are independent copies of a random variable $\xi$ with zero mean, variance one and subexponential tails with exponent $\alpha>0$. Let $N_{n+1}$ be a unit normal vector to $L_n$.
\begin{itemize}
 	\item Then there are constants $c_1,c_2,c_3\in(0,\infty)$ such that
		\begin{equation*}
			\Pro \biggl(\|N_{n+1}\|_\infty \geq \sqrt{\frac{m}{n}}\biggr)\leq c_2\, n^2\exp\Bigl(-c_3\Bigl(\frac{m}{\log n}\Bigr)^{\frac{1}{2\alpha+1}}\Bigr),
		\end{equation*}
		for every $m\geq c_1(\log n)^{2\alpha+2}$.
		As a consequence, with probability $1- c_4 n^{-10}$, say, we have that
		\begin{equation}\label{eq:BerryEsseenProof3}
			\|N_{n+1}\|_\infty \leq  c_5 \frac{(\ln n)^{\alpha+1}}{\sqrt{n}},
		\end{equation}
		for some constant $c_4,c_5>0$.

 	\item If $(V_n)_{n\in\mathbb{N}}$ is a fixed sequence of unit random vectors $V_{n}\in \S^{n-1}$, then
		\begin{equation}\label{eqn:VnNn}
			\sqrt{n}\langle V_{n}, N_{n}\rangle \overset{d}\longrightarrow Z, \quad \text{as $n\to \infty$.}
		\end{equation}
\end{itemize}

\end{lemma}
\begin{rmk}\label{rmk:note}
	Note that \cite[Theorem 1.4]{NguyenVuNormal} is stated for subgaussian random variables. By \cite[Remark 2.3]{NguyenVuNormal} the theorem holds true also for random variables with subexponential tails with exponent $\alpha>0$, but one has to be more generous with the estimates.
\end{rmk}

With probability one the vectors $Y_1,\dotsc,Y_n$ span a random $n$-dimensional linear subspace $L_n$ in $\R^{n+1}$ and we denote by $N_{n+1}\in\S^n$ a unit normal vector to $L_n$. Then
\begin{align*}
	\operatorname{dist}(U_{n+1},L_n) &= |\langle U_{n+1}, N_{n+1}\rangle|
	&\text{and}&&
	\operatorname{dist}(Y_{n+1},L_n) &= |\langle Y_{n+1}, N_{n+1}\rangle|,
\end{align*}
where $\langle{\,\cdot\,},{\,\cdot\,}\rangle$ is the standard scalar product in $\R^{n+1}$.

For the first estimate we set
\begin{equation*}
	V_{n+1} := \frac{1}{\sqrt{n+1}} U_{n+1},
\end{equation*}
use \eqref{eqn:VnNn}, and conclude by the continuous mapping theorem \cite[Lemma 4.3]{Kallenberg}, applied to the absolute-value function, that
\begin{equation}
	\frac{\ln \operatorname{dist}(U_{n+1},L_n)}{\sqrt{\frac{1}{2} \ln n}} \overset{d}\longrightarrow 0
\end{equation}
as $n\to\infty$. This settles the first case of \eqref{eqn:proofAii}.

The second statement of \eqref{eqn:proofAii} also follows by Slutsky's theorem once we show that
\begin{equation}\label{eqn:proofAsecondcase}
	\operatorname{dist}(Y_{n+1},L_n) \overset{d}\longrightarrow |Z|, \quad \text{as $n\to \infty$.}
\end{equation}
To prove this we use Lemma \ref{lem:Berry-Esseen-sup} and the first part of Lemma \ref{le:NguyenVuNormal}. We condition on $L_n$ to fix $N_{n+1}$ and combine \eqref{eq:BerryEsseenProof3} with \eqref{eq:BerryEsseenProof2}. To be more precise,
we have
\begin{align*}
	&\Pro\left(\operatorname{dist}(Y_{n+1},L_n)\leq t\right) 
	= \E \, \Pro\left(\left|\langle Y_{n+1}, N_{n+1} \rangle\right|\leq t \vert L_n) \right)  .
\end{align*}
Here $\Pro(\,{\cdot}\,|L_n)$ denotes the conditional probability given $L_n$, where $L_n$ is the linear space spanned by $Y_1,\dotsc,Y_n\in \R^{n+1}$ and $\E$ denotes expectation with respect to $Y_1, \dots, Y_n$. If we condition on $L_n$, then $N_{n+1}$ is a fixed unit vector in $\R^{n+1}$ and we may apply the Berry-Esseen inequality \eqref{eq:BerryEsseenProof2} with $V_{n+1}=N_{n+1}$ there to deduce that
\begin{equation*}
	\sup_{t\in \R} \, \bigl| \Pro\left(\left|\langle Y_{n+1}, N_{n+1} \rangle\right|\leq t\, |\, L_n \right) - \Pro\left(|Z|\leq t\right)\bigr| \leq c \|N_{n+1}\|_{\infty}.
\end{equation*}
Moreover, from \eqref{eq:BerryEsseenProof3} we conclude that there exist constants $c_1,c_2,c_3\in(0,\infty)$ such that
\begin{equation*}
	0\leq \sup_{t\in \R} \, \bigl| \Pro\left(\left|\langle Y_{n+1}, N_{n+1} \rangle\right|\leq t\, |\, L_n \right) - \Pro\left(|Z|\leq t\right)\bigr| \leq c_2 \frac{(\ln n)^{c_3}}{\sqrt{n}}
\end{equation*}
holds true with probability $1-c_1 n^{-10}$ for sufficiently large $n$.
Hence we have
\begin{align*}
\KS  (\operatorname{dist}(Y_{n+1},L_n),|Z|) 
	& \leq
	\E \, \sup_{t \in \R} \, \left| 
	\Pro\left(\left|\langle Y_{n+1}, N_{n+1} \rangle\right|\leq t \vert L_n\right) 
	- \Pro\left(|Z|\leq t\right) 
	\right|
	\\ & \leq
	c_2\frac{(\ln n)^{c_3}}{\sqrt{n}} (1-c_1 n^{-10}) + 2c_1 n^{-10}.
\end{align*}
Finally, we notice that the last expression tends to zero, as $n\to\infty$. This yields \eqref{eqn:proofAsecondcase} and completes the proof.\hfill $\Box$

% % % % % % % % % % % % % %% % % % % % % % % % % % % %% % % % % % % % % % % % % %% % % % % % % % % % % % % %
% % % % % % % % % % % % % %% % % % % % % % % % % % % %% % % % % % % % % % % % % %% % % % % % % % % % % % % %
% % % % % % % % % % % % % %% % % % % % % % % % % % % %% % % % % % % % % % % % % %% % % % % % % % % % % % % %
% % % % % % % % % % % % % %% % % % % % % % % % % % % %% % % % % % % % % % % % % %% % % % % % % % % % % % % %

\section{Proof of Theorem \ref{thm:CLTlpCASE}}\label{sec:ProofCLTlpCase}

\subsection{The probability measures \texorpdfstring{$\nu=\nu_n(m_n,p)$ on the $\ell_p$-ball}{on the Lp-ball}}\label{sec:nu}

Denote for each $i\in\{1,\dotsc,n\}$ by $G^i_1,\ldots,G^i_n$ random variables with density
\begin{align}\label{eqn:const_a}
	t\mapsto\frac{e^{-|t|^p/a}}{2a^{1/p}\Gamma(1+1/p)},\qquad
	\text{where}\qquad 
	a 
	  = \left(\frac{\Gamma(1/p)}{\Gamma(3/p)}\right)^{p/2},
	\qquad t\in\R.
\end{align}
Note that $G^i_j$ has zero mean and variance one and subexponential tails with exponent $\alpha=1/p$.
In addition, let $m_n\in[0,\infty)$ and, for each $i\in\{1,\ldots,n\}$, $Q_{i}$ be random variables which are gamma distributed with shape $m_n$ and rate $1/a$. More specifically this means that $Q_i$ has density $t\mapsto a^{-m_n}\Gamma(m_n)^{-1}t^{m_n-1}e^{-t/a}$ for $t>0$, provided that $m_n > 0$, and we use the convention that $Q_i=0$ with probability one in case that $m_n=0$. We shall assume that all the random variables we are considering are independent.

Next, we define the random vectors $X_1,\ldots,X_n\in\R^n$ by putting
\begin{equation*}
	X_i := \frac{G_i}{(\|G_i\|_p^p+Q_{i})^{1/p}}, \qquad\text{where} \qquad G_i:=(G^i_1,\ldots,G^i_n),
\end{equation*}
for $i\in\{1,\ldots,n\}$. By $\nu_n(m_n,p)$ we denote the distribution of the random variables $X_i$ on the $\ell_p$-ball $B_p^n$. Finally, we let $\Xi_n(K,\nu)$ be the random convex body that is generated by a fixed convex body $K\subset \R^n$ and the independent random points $X_1,\ldots,X_n$ for a distribution $\nu=\nu_n(m_n,p)$, which we consider to be fixed in this section.

\subsection{Proof of Theorem \ref{thm:CLTlpCASE}}

Recall that by \eqref{eqn:randbody_vol} we have
\begin{align}\label{eqn:CLTlpFirstStep}
\begin{split}
	\vol_n(\Xi_n(K_n,\nu)) &=  \vol_n(K_n) \left|\det(X_1|\dotsc|X_n)\right|\\
	&= \vol_n(K_n)  \left|\det(G_1|\dotsc|G_n)\right|\prod_{i=1}^n \big(\|G_i\|_p^p+Q_i\big)^{-1/p},
\end{split}
\end{align}
which yields
\begin{equation*}
	\ln \vol_n(\Xi_n(K_n,\nu)) = \ln\,\left|\det(G_1|\dotsc|G_n)\right| + \ln \vol_n(K_n) - \frac{1}{p} \sum_{i=1}^n \ln\left(\|G_i\|_p^p+Q_i\right).
\end{equation*}
Now we may apply the central limit theorem for the log-determinant Theorem \ref{thm:CLTDeterminantNguyenVu} and obtain
\begin{equation}\label{eq:lpCaseCLTDeterminant}
	O_n:=\frac{\ln\,\left|\det(G_1|\dotsc|G_n)\right| - \frac{1}{2}\ln\, (n-1)!}{\sqrt{\frac{1}{2}\ln n}} \overset{d}\longrightarrow Z,\quad \text{as $n\to\infty$},
\end{equation}
and moreover $\KS(O_n,Z) \leq (\ln n)^{-\frac{1}{3}+o(1)}$.
To complete the proof of Theorem \ref{thm:CLTlpCASE} we need to show that
\begin{equation}\label{eqn:CLTlpSecondStep}
	P_n:=\frac{\frac{1}{p} \sum\limits_{i=1}^n \ln\left(\|G_i\|_p^p+Q_i\right) - \frac{n}{p} \ln(a(m_n+\frac{n}{p}))}{\sqrt{\frac{1}{2}\ln n}} \overset{d}\longrightarrow 0,
\end{equation}
and then apply Slutsky's theorem. Indeed, putting together \eqref{eqn:CLTlpFirstStep} with \eqref{eq:lpCaseCLTDeterminant} and \eqref{eqn:CLTlpSecondStep} implies that
\begin{align*}
	&\frac{\ln \vol_n(\Xi_n(K_n,\nu)) - \ln \vol_n(K_n) - \frac{1}{2} \ln\, (n-1)! + \frac{n}{p}\ln(a(m_n+\frac{n}{p}))}{\sqrt{\frac{1}{2}\ln n}}\\
	&\quad = 
	O_n-P_n
	 \overset{d}\longrightarrow Z,\quad \text{as $n\to \infty$,}
\end{align*}
as desired.
For all $\varepsilon>0$ we have that
\begin{equation*}
	\KS(O_n-P_n,Z) \leq \KS(O_n,Z) + \Pro(|P_n|>\varepsilon) + \varepsilon,
\end{equation*}
see for example \cite[Lemma 4.1]{APT17CLT}.
Hence, once we show that
\begin{equation}\label{eqn:CLTlpSecondStep2}
	\Pro(|P_n|>(\ln n)^{-\frac{1}{3}}) \leq (\ln n)^{-\frac{1}{3}+o(1)},
\end{equation}
we may conclude, by setting $\varepsilon=(\ln n)^{-\frac{1}{3}}$ and applying Theorem \ref{thm:CLTDeterminantNguyenVu}, that
\begin{equation*}
	\KS(O_n-P_n,Z) \leq (\ln n)^{-\frac{1}{3}+o(1)}.
\end{equation*}
This will finish the proof of Theorem \ref{thm:CLTlpCASE}.

\subsection{Proof of \texorpdfstring{\eqref{eqn:CLTlpSecondStep2}}{(\ref{eqn:CLTlpSecondStep2})}}

We observe that the representation
\begin{equation*}
	\|G_i\|_p^p = \sum_{j=1}^n|G^i_j|^p
\end{equation*}
and the semigroup property of the gamma distributions imply that $\|G_i\|_p^p$ is gamma distributed with shape $n/p$ and rate $1/a$, i.e., the Lebesgue density of $\|G_i\|_p^p$ on $\R$ is given by
\begin{equation*}
	\frac{1}{a^{\frac{n}{p}}\Gamma(\frac{n}{p})}x^{\frac{n}{p}-1}e^{-\frac{x}{a}},\qquad x>0.
\end{equation*}
Recalling that by assumption $Q_i$ is gamma distributed with shape $m_n$ and rate $1/a$ we find that $\|G_i\|_p^p+Q_i$ is also gamma distributed with parameter shape $m_n+\frac{n}{p}$ and rate $1/a$. Hence, $\ln(\|G_i\|_p^p+Q_i)$ is log-gamma distributed with Lebesgue density on $\R$ given by
\begin{equation*}
	\frac{1}{a^{m_n+\frac{n}{p}}\Gamma(m_n+\frac{n}{p})} e^{x(m_n+\frac{n}{p})-\frac{e^x}{a}},\qquad x>0.
\end{equation*}
In particular, by direct computation, we find that
\begin{equation}\label{eq:Mun}
	\mu_n:=\E \ln(\|G_i\|_p^p+Q_i) = \psi\left(m_n+\frac{n}{p}\right)+\ln a,% \sim \ln\left(a\left(m_n+\frac{n}{p}\right)\right),
\end{equation}
where $\psi(x):=\frac{d}{dx}\ln\Gamma(x)$ is the digamma function (see e.g.\ \cite[page 259]{AbramowitzStegun}). 
Similarly, for the variance, one has that
\begin{equation*}
	\sigma_n^2:=\Var\ln(\|G_i\|_p^p+Q_i) = \psi_1\left(m_n+\frac{n}{p}\right),% \sim \frac{1}{m_n+\frac{n}{p}},
\end{equation*}
with $\psi_1(x):=\frac{d}{dx}\psi(x)$ being the trigamma function (see e.g.\ \cite[page 260]{AbramowitzStegun}).

This implies that the auxiliary random variables
\begin{equation*}
	A_n:=\left[\frac{1}{p}\sum_{i=1}^n\ln \left(\|G_i\|_p^p+Q_i\right)\right] - \frac{n}{p}\mu_n
\end{equation*}
satisfy
\begin{equation*}
	\E A_n = 0 \qquad \text{and} \qquad \Var A_n = \frac{n}{p^2} \psi_1\left(m_n+\frac{n}{p}\right).
\end{equation*}
The asymptotic expansions of the digamma and trigamma functions are 
\begin{align*}
	\psi(x) &= \ln x - \frac{1}{2x} - \frac{1}{12x^2} + o(x^{-2}),&
	\psi_1(x) &= \frac{1}{x} + \frac{1}{2x^2} + o(x^{-2}),
\end{align*}
for $x\to \infty$ (see \cite[page 260]{AbramowitzStegun}).
Hence
\begin{align*}
	\frac{\frac{n}{p}\mu_n - \frac{n}{p} \ln(a(m_n+\frac{n}{p}))}{\sqrt{\frac{1}{2} \ln n}} 
	= \frac{1}{\sqrt{2}} \left(1+p\frac{m_n}{n}\right)^{-1} (\ln n)^{-\frac{1}{2}} + o((\ln n)^{-\frac{1}{2}})\longrightarrow 0,
\end{align*}
for $n\to \infty$, and
\begin{equation*}
	\Var A_n = \frac{n}{p} \frac{1}{pm_n+n} + \frac{n}{2} \frac{1}{(pm_n+n)^2} + o(n^{-1}).
\end{equation*}
In particular, for all choices of $m_n$ we find that $\Var A_n = O(1)$.
By the triangle inequality we have
\begin{align*}
	|P_n| \leq \frac{|A_n|}{\sqrt{\frac{1}{2}\ln n}} + \frac{1}{\sqrt{2}} \left(1+p\frac{m_n}{n}\right)^{-1} (\ln n)^{-\frac{1}{2}} + o((\ln n)^{-\frac{1}{2}})
\end{align*}
and therefore there exists $c_1>0$ such that
\begin{align*}
	|P_n|\leq \frac{|A_n|+c_1}{\sqrt{\frac{1}{2}\ln n}},
\end{align*}
for all $n$ large enough.
By the Chebyshev inequality this yields
\begin{align*}
	\Pro(|P_n|>(\ln n)^{-\frac{1}{3}}) 
	&\leq \Pro\left(|A_n| > \frac{1}{\sqrt{2}} (\ln n)^{\frac{1}{6}} - c_1\right) \\
	&\leq \Pro\left(|A_n| > \frac{1}{2}(\ln n)^{\frac{1}{6}}\right)
	\leq \frac{\sqrt{2}\Var A_n}{(\ln n)^{\frac{1}{3}}} = (\ln n)^{-\frac{1}{3}+o(1)},
\end{align*}
for all $n$ large enough.
Thus, \eqref{eqn:CLTlpSecondStep2} holds true and the proof is complete. \hfill $\Box$

\begin{rmk}
	More general distributions for the random variables $Q_i$ are possible. For example, our proof shows that, 
	as long as $Q_i$ is a non-negative random variable and ${\rm Var} \ln(\|G_i\|_p^p+Q_{i}) = O\bigl({\rm Var} \ln(\|G_i\|_p^p)\bigr) = O(p/n)$, 
	we have a CLT as above with the same scaling factor and a suitably modified final centering term.
\end{rmk}

% % % % % % % % % % % % % %% % % % % % % % % % % % % %% % % % % % % % % % % % % %% % % % % % % % % % % % % %
% % % % % % % % % % % % % %% % % % % % % % % % % % % %% % % % % % % % % % % % % %% % % % % % % % % % % % % %

\newpage

\end{document}